\documentclass[11pt,twoside]{article}
\usepackage{mathbbold}
\usepackage{amsmath, amssymb, mathrsfs, graphicx}

\def\scr{\mathscr}

\def\az{\alpha}  \def\bz{\beta}
    \def\dz{\delta}
\def\ez{\eta}    \def\fz{\varphi}
\def\gz{\gamma}  
 
\def\nz{\nu}

        \def\sz{\sigma}
        
\def\vz{\varepsilon}

\def\qd{\quad}
\def\qqd{\qquad}

\newcommand{\mathsym}[1]{{}}

\def\le{\leqslant}
\def\ge{\geqslant}

\font\cms=cmss9 scaled \magstep1

\def\nnd{\noindent}

\def\thm{\nnd\bg{thm1}}
\def\crl{\nnd\bg{crl1}}
\def\lmm{\nnd\bg{lmm1}}
\def\prp{\nnd\bg{prp1}}

\def\xmp{\nnd\bg{xmp1}}

\def\rmk{\nnd\bg{rmk1}}

\def\dethm{\end{thm1}}
\def\decrl{\end{crl1}}
\def\delmm{\end{lmm1}}
\def\deprp{\end{prp1}}
\def\dexmp{\end{xmp1}}

\def\dermk{\end{rmk1}}

\def\prf{\medskip \noindent {\bf Proof}. }
\def\qed{\text{\quad $\square$}}
\def\deprf{\qed\medskip}
\def\bg{\begin}
\def\be{\bg{equation}}
\def\de{\end{equation}}

\def\dear{\end{eqnarray}}
\def\lb{\label}
\def\ct{\cite}

\newcommand{\rf}[2]{[\ref{#1}; #2]}

\def\den{\end{enumerate}}

\def\d{\text{\rm d}}

\allowdisplaybreaks[4]

\begin{document}

\thispagestyle{empty}
\renewcommand{\thefootnote}{\fnsymbol{footnote}}


\vspace*{.5in}
\begin{center}
{\bf\Large The optimal constant in Hardy-type inequalities}
\vskip.15in {Mu-Fa Chen}
\end{center}
\begin{center} (Beijing Normal University, Beijing 100875, China)\\
\vskip.1in March 19, 2013
\end{center}
\vskip.1in

\markboth{\sc Mu-Fa Chen}{\sc Hardy-type inequalities}



\date{}




\begin{abstract}To
estimate the optimal constant in Hardy-type inequalities, some variational formulas and approximating procedures are introduced. The known basic estimates
are improved considerably. The results are illustrated by typical examples. It is shown that the sharp factor is meaningful for each finite interval and a classical sharp model is re-examined.
    \end{abstract}

\section{Introduction}

For given two Borel measures $\mu$ and $\nu$ on an interval $(-M, N)\,(M, N\le \infty)$,
the Hardy-type inequality says that the $L^q(\mu)$-norm of each
absolutely continuous function $f$ is controlled from above
by the $L^p(\nu)$-norm of its derivative $f'$ up to a constant factor $A$:
\begin{align}
&\bigg(\int_{-M}^N |f|^q\d\mu\bigg)^{1/q}\!\!\le A
 \bigg(\int_{-M}^N |f'|^p\d{\nu}\bigg)^{1/p},\lb{01}\\
 &\text{i.e.,}\; \|f\|_{\mu, q}\le A \|f'\|_{\nu, p},\qqd p, q\in (1, \infty).\nonumber
 \end{align}
The inequalities have been well studied in the past decades, cf. \ct{ok-90, cmp-07, cmf13}. In particular, the following basic estimates for the optimal constant $A$ in (\ref{01}) are well known:
\be B\le A \le k_{q,p} B\lb{02},\de
where $B$ is a quantity described by $M, N, \mu, \nu, p$ and $q$, and $k_{q,p}\in [1, 2]$
 is a constant factor (cf. (\ref{03}), (\ref{22}) and (\ref{23}) below. See also Appendix for
 more details). The goal of this paper
is to show that there is still a room for improvements of (\ref{02}). Such a qualitative study is valuable since the optimal constant $A$ describes the speed of some type of stability
(cf. \rf{cmf05a}{Chapter 6} and references therein). We begin our story with
an example which is typical in the sense that it is the only one,
except the special case that $p=q=2$, we have known so far
for having the exact constant $A$ for all $p, q\in (1, \infty)$
(see also Example \ref{t2-5} and Proposition \ref{ta-5} for additional information).

There are mainly four different types of boundary conditions in (\ref{01}). We concentrate in the paper only on the one vanishing at $-M$. The other cases will be handled subsequently. In particular, the results of this paper
are extended by \ct{liao-14} to the discrete context.

\xmp\lb{t1-1}{\rm Let $(-M, N)=(0, 1)$ and $\d\mu=\d\nu=\d x$. Then

(1) {\it the optimal constant $A$} in (\ref{01}) is
$$
A=\frac{p^{\frac{1}{q}} q^{1-\frac{1}{p}} (pq+p-q)^{\frac{1}{p}-\frac{1}{q}}
}
{(p-1)^{\frac 1 p} \text{\rm B}\!\left(\frac{1}{q},\,1-\frac{1}{p}\right)
 },\qqd p, q\in (1, \infty),$$
 where $\text{\rm B}(\az, \bz)$ is the Beta function
$$\text{\rm B}(\az, \bz)\!=\!\int_0^1\!\! s^{\az-1}(1-s)^{\bz-1}\d s=\frac {1}{\az} \int_0^1 \!\! \big(1-t^{\frac 1 {\az}}\big)^{\bz-1}\d t\qd\text{(change variable $t=s^{\az}$)}.$$
In particular, if $q=p$, then
 $$A=\frac{p}{\pi (p-1)^{1/p}}\sin\frac \pi p.$$
More particularly, $A=2/\pi$ if $q=p=2$.

(2) {\it Basic estimates}. The constants used for the basic estimates in (\ref{02}) for $q\ge p$ are as follows:
\begin{align}
&\text{\hspace{-2em}}B=   \frac{p^{\frac{1}{q}}\big((p-1)q\big)^{1-\frac 1 p}}
{(pq+p-q)^{1-\frac 1 p + \frac 1 q}},\nonumber\\
&\text{\hspace{-2em}}k_{q, p}\!=\!\!\Bigg[\frac{\Gamma\big(\frac{pq}{q-p}\big)}{\Gamma\big(\frac{q}{q-p}\big)\Gamma\big(\frac{p(q-1)}{q-p}\big)}\Bigg]^{1/p-1/q}
\!\!=\!\bigg[\frac{q-p}{p \text{\rm B}\!\Big(\frac{p}{q-p},\, \frac{p(q-1)}{q-p}\Big)}\Bigg]^{1/p-1/q}\!\!\!\!,\qquad
q\!>\!p, \lb{03}\end{align}
cf. \rf{cmf13}{Example 1.12 and references therein}. In particular, if $q=p$, then
$$B=\big(1/ p\big)^{1/p}
\big(1/ {p^*}\big)^{1/p^*},\qqd k_{p, p}=p^{1/p}{p^*}^{1/p^*},$$
where $p^*$ is the conjugate of $p$: $1/p+1/p^*=1$.

(3) {\it Improvements}. As the first step of our approximating procedures introduced in the paper, we have
new upper and lower bounds $\dz_1$ and ${\bar\dz}_1$, respectively. Besides, we also have
another upper estimate $A^*$. More precisely, the basic estimates in (\ref{02}) are improved in the paper by
\be B\le {\bar\dz}_1\le A \le A^*\le \dz_1\le k_{q, p}B, \lb{04}\de
where
\begin{align}
\mbox{\hskip-0.5em}{\bar\dz}_1\!&=\!\frac{p^{1/q}((p-1)(q+1))^{1-1/p}}{(pq+p-q)^{1-1/p+1/q}},\\
\mbox{\hskip-0.5em}A^*\!&=\bigg(\frac{p^*}{q}\bigg)^{1/ q}
\bigg(\frac{p^*+q}{\pi p^*}\sin \frac{\pi p^*}{p^*+q}\bigg)^{{1}/{p^*}+ 1/q},\\
\mbox{\hskip-0.5em}{\dz_1}\!&=\!\frac{1}{({q\gz^*}/{p^*}\!+\!1)^{1/ q}}
\bigg\{\!\sup_{x\in (0, 1)}\frac{1}{x^{\gz^*}}\!\!\int_0^x\!\!\!\big(1\!-\!y^{{q\gz^*}/{p^*}+1}\big)^{{p^*}/{q}}\d y\bigg\}^{{1}/{p^*}}\!\!\!,\;\;
\gz^*\!\!:=\!\frac{q}{p^*\!+\!q}.
\end{align}
In the last formula, the function under $\sup_{x\in (0, 1)}$ is unimodal on $(0, 1)$, its
integral term is indeed an incomplete Beta function:
$$\text{B}(x, \az, \bz)=\int_0^x s^{\az-1}(1-s)^{\bz-1}\d s
=\frac{1}{\az}\int_0^{x^{\az}} \big(1-t^{1/\az}\big)^{\bz-1}\d t.$$
Note that ${\bar\dz}_1$ is very much the same as $B$: the factor $q+1$ in ${\bar\dz}_1$
is replaced by $q$ in $B$. Besides, $A^*=A$ if $q=p$. Except these facts, the comparison of the quantities in (\ref{04})
are non-trivial, as shown by Figures 1--4.}
\dexmp

(4) {\it Figures}. First, consider the case that $q=p$. Figure 1 shows the basic estimates
of the optimal constant $A$.
\begin{center}{\includegraphics[width=11.0cm,height=7.5cm]{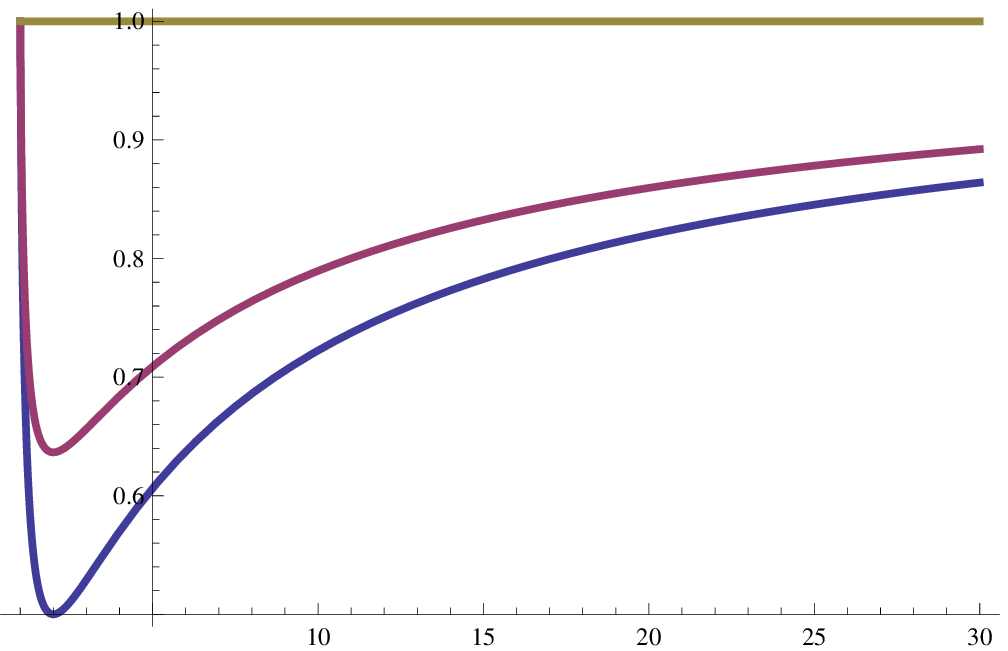}}\end{center}
\nnd{\bf Figure 1}\quad The middle curve is the exact value of $A$.
The top straight line and the bottom curve consist of the basic estimates of $A$.

Our improved upper bound ${\dz}_1$ and lower one ${\bar\dz}_1$ are added to Figure 1, as shown in Figure 2.
\begin{center}{\includegraphics[width=11.0cm,height=7.5cm]{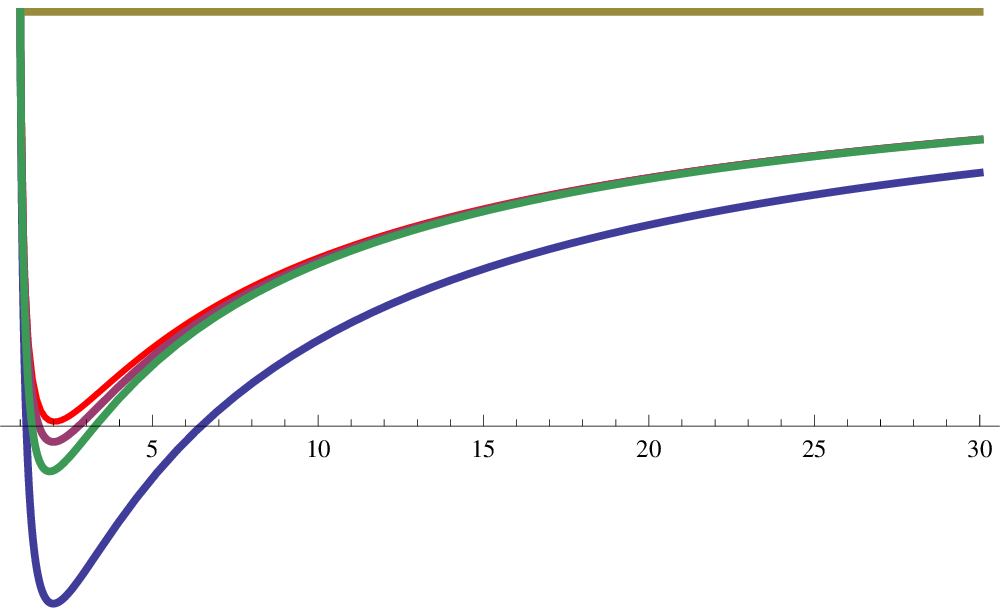}}\end{center}
\nnd{\bf Figure 2}\quad The new bounds ${\dz}_1$ and ${\bar\dz}_1$ are almost overlapped with the exact
value $A$ except in a small neighborhood of $p=2$, $\dz_1$ is a little bigger,
and ${\bar\dz}_1$ is a little smaller than $A$.

Next, consider the case that $q>p$. For convenience, we rewrite $q$ as $p+r$, where
$r$ varies over $(0, 15)$. The six quantities in (\ref{04})
are shown in Figures 3 and 4 according to $p=2$ and $p=5$, respectively.
\begin{center}{\includegraphics[width=11.0cm,height=7.5cm]{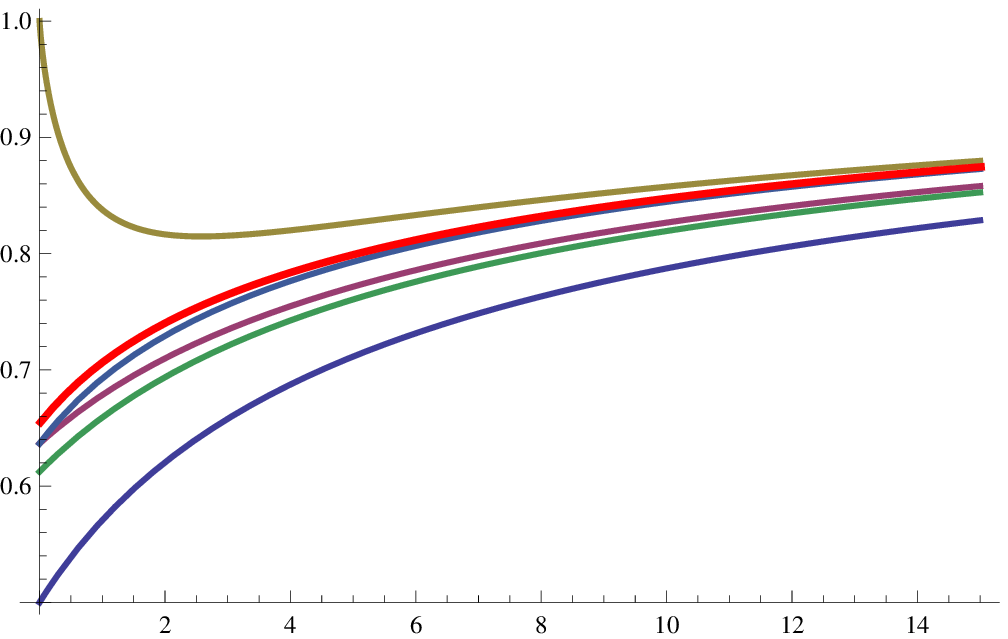}}\end{center}
\nnd{\bf Figure 3}\quad The six curves from top to bottom are $k_{q,p}B$,
$\dz_1$, $A^*$, $A$ (the exact value), ${\bar\dz}_1$, and $B$, respectively,
for $p=2$, $q=p+r$, and $r\in (0, 15)$.

In view of Figures 1 and 2, it is clear that the six curves should be closer for larger $p$.
\begin{center}{\includegraphics[width=11.0cm,height=7.5cm]{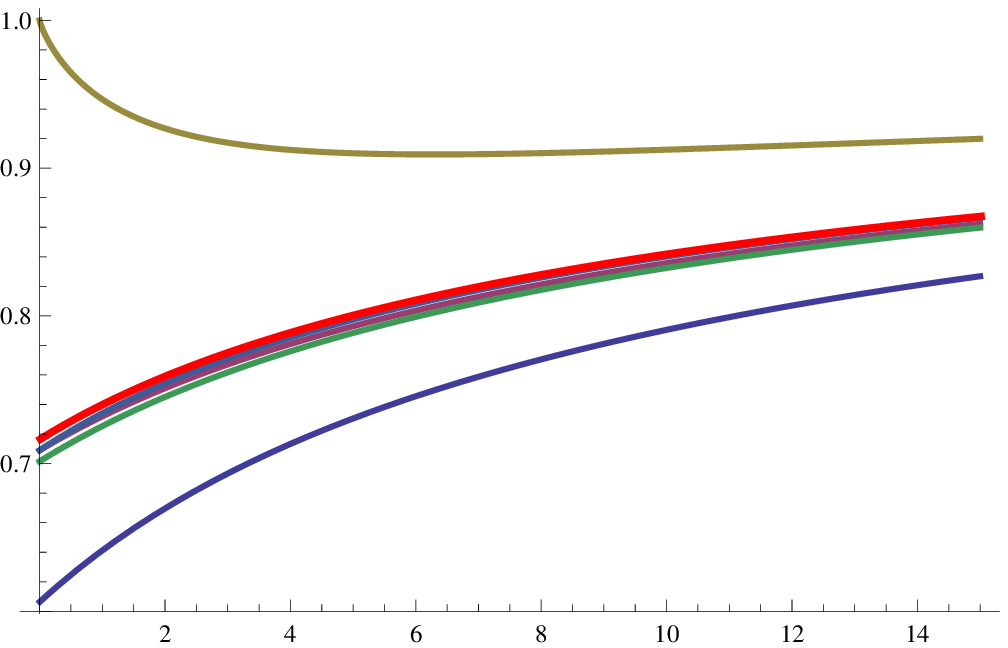}}\end{center}
\nnd{\bf Figure 4}\quad The only change of this figure from the last one is
replacing $p=2$ by $p=5$. Certainly, the six curves are located in the same order.
Except the basic estimates, the other four curves are almost overlapped.

Figures 1--4 illustrate the effectiveness of our improvements. It is surprising and
unexpected that the new estimates can be so closed to the exact value. The general results are
presented in the next section. Their proofs are given in Section 3.
In Appendix, we will come back to study the basic estimates (\ref{02}) and the optimal factor $k_{q, p}$.

To conclude this section, we make some historical remarks on Example \ref{t1-1}.
The optimal constant given in the example was presented in \rf{tg76}{page 357} with optimizer but without details.
The detailed proofs were presented in \ct{dm99} and \ct{aprs}.
We mention that the boundary condition used in the cited papers are vanishing at both
endpoints. This is using the Dirichlet boundaries at two endpoints. The result is the same
if we replace Dirichlet boundaries with Neumann ones (cf. \ct{cmf13}). However, as mentioned before, we consider
only the Dirichlet boundary at the left-endpoint in this paper.
Then the optimal constant here is a double of those given in the cited papers.
For more recent progress on $p$-Laplacian (which is an alternative description of
the Hardy-type inequality in the case of $q=p$), one may refer to the book \ct{LE-11}
and references therein. Actually, in this case, the story is now quite complete.
The new progress will be published elsewhere in \ct{cwz-13}.

\section{Main results}

From now on, for simplicity, we fix $(-M, N)=(0, D)$, $D\le +\infty$. Set
$$ \aligned
&{\scr A}[0, D]=\{f: f\text{ is absolutely continuous in }[0, D]\,(\text{or }[0, D)\text{ if }D=\infty)\},\\
&{\scr A}_0[0, D]=\{f\in {\scr A}[0, D]: f(0)=0\}.
\endaligned$$
Then the optimal constant $A$ in (\ref{01}) is described by the following
classical variational formula
\be A=\sup_{f\in {\scr A}_0[0, D],\; \|f'\|_{\nu, p}\in (0, \infty)}
\frac{\|f\|_{\mu, q}}{\|f'\|_{\nu, p}}\lb{08}\de

To state our results, we need some notation.
Denote by $v$ the density of the absolutely continuous part of $\nu$
with respect to the Lebesgue measure $\d x$ and let
$${\hat v}= v^{-\frac{1}{p-1}}=v^{1-p^*}.$$
For upper estimates, define two operators $I\!I^*$ and $I^*$:
\begin{align}
I\!I^*(f)(x)&=\frac 1 {f(x)}\int_0^x \d y \,{\hat v}(y) \bigg(\int_y^D f^{q/p^*}\d \mu\bigg)^{{p^*}/{q}},\qqd x\in (0, D),\\
I^*(f)(x)&=\frac{\hat v}{f'}(x)
\bigg(\int_x^D f^{q/p^*}\d \mu\bigg)^{p^*/q},\qqd x\in (0, D),
\end{align}
with domains
\begin{align}
{\scr F}_{I\!I}&=\{f: f(0)=0, f>0 \text{ on } (0, D)\},\\
{\scr F}_{I}&=\{f: f(0)=0, f'>0 \text{ on } (0, D)\},
\end{align}
respectively.
For lower estimates, we need different operators:
\begin{align}
I\!I(f)(x)&=\frac{1}{f(x)}\int_0^x \d y\,{\hat v}(y)\bigg(\int_y^D f^{q-1}\d\mu\bigg)^{p^*-1},
\qqd x\in (0, D),\\
I(f)(x)&=\frac{{\hat v}}{f'}(x)\bigg(\int_x^D f^{q-1}\d\mu\bigg)^{p^*-1},\qqd x\in (0, D).
\end{align}
When $q=p$, we have $I\!I=I\!I^*$ and $I=I^*$.
To avoid the non-integrability problem, the domain of $I\!I$ and $I$ have to be modified
from ${\scr F}_{I\!I}$ and ${\scr F}_{I}$:
\begin{align}
{\widetilde{\scr F}}_{I\!I}&=\{f\in {\scr F}_{I\!I}: \exists x_0\in (0, D]\text{ such that }
        f=f(\cdot\wedge x_0)\text{ and moreover } \nonumber\\
&\text{\hskip 6em} f I\!I(f)\in L^q(\mu) \text{ if }x_0=D\},\\
{\widetilde{\scr F}}_{I}&=\{f: f(0)=0, \; \exists x_0\in (0, D]\text{ such that }
        f=f(\cdot\wedge x_0),\; f'>0\text{ on $(0, x_0)$,} \nonumber\\
&\text{\hskip 5.5em}\; \text{and moreover } f I\!I(f)\in L^q(\mu) \text{ if }x_0=D\},
\end{align}
where $\az\wedge\bz=\min\{\az, \bz\}$ and similarly, $\az\vee\bz=\max\{\az, \bz\}$.
Thus, the operators we actually use for the lower estimates are modified from $I\!I$
and $I$ as follows: when $f=f(\cdot\wedge x_0)$, set
\begin{align}
I\!I\, \tilde{} (f)(x)&\!=\!I\!I(f)(x\wedge x_0)=\frac{1}{f(x)}\int_0^{x\wedge x_0} \d y\,{\hat v}(y)\bigg(\int_y^D f^{q-1}\d\mu\bigg)^{p^*-1},\\
I\, \tilde{}(f)(x)&\!=\!I(f)(x\wedge x_0)\!=\!\frac{{\hat v}}{f'}(x\wedge x_0)\bigg(\int_{x\wedge x_0}^D f^{q-1}\d\mu\bigg)^{p^*-1}\!\!\!,\qd x\in (0, D).
\end{align}
Here we adopt the usual convention that $1/0=\infty$.

We can now state our variational result.

\thm\lb{t2-1}{\cms For the optimal constant $A$ in (\ref{01}), we have
\bg{itemize}
\item[(1)] upper estimate:
\be A\le \inf_{f\in {\scr F}_{I\!I}}\Big[\sup_{x\in (0, D)}I\!I^*(f)(x)\Big]^{1/p^*}= \inf_{f\in {\scr F}_{I}}\Big[ \sup_{x\in (0, D)}I^*(f)(x)\Big]^{1/p^*}\lb{19} \de
for $q\ge p$, and
\item[(2)] lower estimate:
\be A\ge  \sup_{f\in {\widetilde{\scr F}}_{I\!I}} \|f I\!I\, \tilde{}(f)\|_{\mu, q}^{1-q/p}\Big(\inf_{x\in (0, D)} I\!I\, \tilde{}(f)(x)\Big)^{(q-1)/p}\lb{20}\de
for $p, q\in (1, \infty)$. In particular, when $q= p$, we have additionally that
\be \sup_{f\in {\widetilde{\scr F}}_{I\!I}}\Big(\inf_{x\in (0, D)} I\!I\, \tilde{}(f)(x)\Big)^{1/p^*}
 =\sup_{f\in {\widetilde{\scr F}}_{I}}\Big(\inf_{x\in (0, D)}I\, \tilde{}(f)(x)\Big)^{1/p^*}.\lb{21}\de
\end{itemize}
}
\dethm

Recall that for general $q\ge p$, $p, q\in (1, \infty)$, the basic estimates read as follows
\be B\le A\le k_{q, p} B  \lb{22}\de
where $k_{q, p}$ is given in (\ref{03}) and
\be B=\sup_{x\in (0, D)} {\hat\nz}(0, x)^{1/p^*} \mu(x, D)^{1/q},\lb{23}\de
here ${\hat\nz}(\az, \bz)\!=\!\int_{\az}^{\bz} {\hat v}$ as usual, and similar for $\mu(\az, \bz)$ (cf. \rf{cmp-07}{pages 45--47}
and Appendix below).

As an application of Theorem \ref{t2-1}, we have the following approximating procedures.

\thm\lb{t2-2}{\cms
\bg{itemize}
\item[(1)] Let $q\ge p$, $p, q\in (1, \infty)$,
\be f_1(x)={\hat\nu}(0, x)^{\gz^*}, \qqd \gz^*=\frac{q}{p^*+q},\de
and define
\bg{align}
&f_{n+1}(x)=f_n I\!I(f_n)=\int_0^x \d y \,{\hat v}(y) \bigg(\int_y^D f_n^{q/p^*}\d \mu\bigg)^{{p^*}/{q}},\qqd n\ge 1,\\
&\dz_n={
\begin{cases}
\Big(\sup_{x\in (0, D)}\frac{f_{n+1}(x)}{f_n(x)}\Big)^{1/p^*},\qd n=1 \text{\cms\; or $n\ge 2$ but $\dz_1<\infty$}\\
\infty, \text{\cms\; $n\ge 2$ and $\dz_1=\infty$}.
\end{cases}
}
\end{align}
Then we have $A\le \dz_n$ for all $n\ge 1$ and $\{\dz_n\}$ is decreasing in $n$.
\item[(2)] Let $p, q\in (1, \infty)$,
\be f_1^{(x_0)}(x)={\hat\nu} (0, x\wedge x_0),\qqd f_{n+1}^{(x_0)}=f_n^{(x_0)} I\!I\, \tilde{}\big(f_n^{(x_0)}\big),\qd n\ge 1\de
and define
\begin{align}
{\tilde \dz}_n&\!\!=\!\!\sup_{x_0\in (0, D]}\big\|f_n^{(x_0)} I\!I\, \tilde{}\big(f_n^{(x_0)}\big)\big\|_{\mu, q}^{1-q/p}\Big(\inf_{x\in (0, D)} I\!I\, \tilde{}\big(f_n^{(x_0)}\big)(x)\Big)^{(q-1)/p},\\
{\bar\dz}_n&\!\!=\!\!\sup_{x_0\in (0, D]} \frac{\big\|f_n^{(x_0)}\big\|_{\mu, q}}{{\big\|f_n'}^{(x_0)}\big\|_{v, p}},\qqd n\ge 1.
\end{align}
Then we have $A\ge {\tilde \dz}_n \vee {\bar\dz}_n$ for all $n\ge 1$.
\end{itemize}
}
\dethm

Actually, in view of Corollary \ref{t2-3} below, we have $\dz_1<\infty$ iff $B< \infty$. When $q=p$, it is known from \ct{cwz-13} that $\{{\tilde \dz}_n\}_{n\ge 1}$ is increasing
in $n$ and ${\tilde \dz}_1\ge B$.

We can now summarize the first step of our approximating procedures
as follows.

\crl\lb{t2-3}{\cms
For $q\ge p>1$, we have
\be B\le {\bar\dz}_1\vee {\tilde\dz}_1\le A\le \dz_1\wedge \big({k}_{q,p} B\big)\le
\dz_1\le {\tilde k}_{q,p} B,\lb{30}\de
where
$$ {\tilde k}_{q,p}= \bigg(1+\frac{q}{p^*}\bigg)^{1/q}
 \bigg(1+\frac{p^*}{q}\bigg)^{1/p^*}\qd\big(\ge k_{q,p}\text{\cms\; if }q\ge p\big).$$
More precisely, let $\fz(x)={\hat\nz}(0, x)$. Then we have
\begin{gather}
\dz_1\!=\!\bigg\{\sup_{x\in (0, D)}
\frac{1}{\fz(x)^{\gz^*}} \int_0^x \d y\,{\hat v}(y) \bigg(\int_y^D \fz^{q \gz^*/p^*}\d \mu\bigg)^{{p^*}/{q}}\bigg\}^{1/p^*}\!\!\!\le\! {\tilde k}_{q, p}B,\lb{31}
\end{gather}
where $\gz^*=\frac{q}{p^*+q}$. Next, we have
\begin{gather}
{\bar\dz}_1=\bigg\{\sup_{x\in (0, D)}\bigg[\frac{1}{\fz(x)^{q/p}}\int_0^x \fz^q\d\mu +\fz(x)^{q/p^*} \mu(x, D)\bigg]\bigg\}^{1/q}\ge B,\lb{32}\\
{\tilde \dz}_1=\sup_{x_0\in (0, D)}
\big\|f_2^{(x_0)}\big\|_{\mu, q}^{1-q/p}\bigg[\frac{f_2^{(x_0)}(x_0)}{\fz(x_0)}\bigg]^{(q-1)/p},
\end{gather}
where
\begin{gather}
f_2^{(x_0)}(x)\!=\!
\!\int_0^{x\wedge x_0}\! \d y\,{\hat v}(y)\bigg[\int_y^{x_0}\!\! \fz^{q-1}\d\mu
   \!+\! \fz(x_0)^{q-1}\mu (x_0, D)\bigg]^{p^*-1}\!\!,\qd x\in [0, D].\nonumber
\end{gather}
}
\decrl

It is known that ${\tilde k}_{p, p}=\lim_{q\downarrow p} k_{q, p}$. When
$q>p$, we have ${\tilde k}_{q, p}>k_{q, p}$. Their small differences are shown by Figure 5.
Their ratios have similar shape as in Figure 5 and are located in $[1,\, 1.23)$ with maximum
$1.2274$ at $(p, q)\approx (2, \, 2+2.5758)$. Besides,
$$\sup_{q\ge p>1}{\tilde k}_{q, p}=\sup_{q\ge p>1}{k}_{q, p}=\sup_{p>1}{\tilde k}_{p, p}
={\tilde k}_{2, 2}=2.$$

{\begin{center}{\includegraphics[width=11.0cm,height=7.5cm]{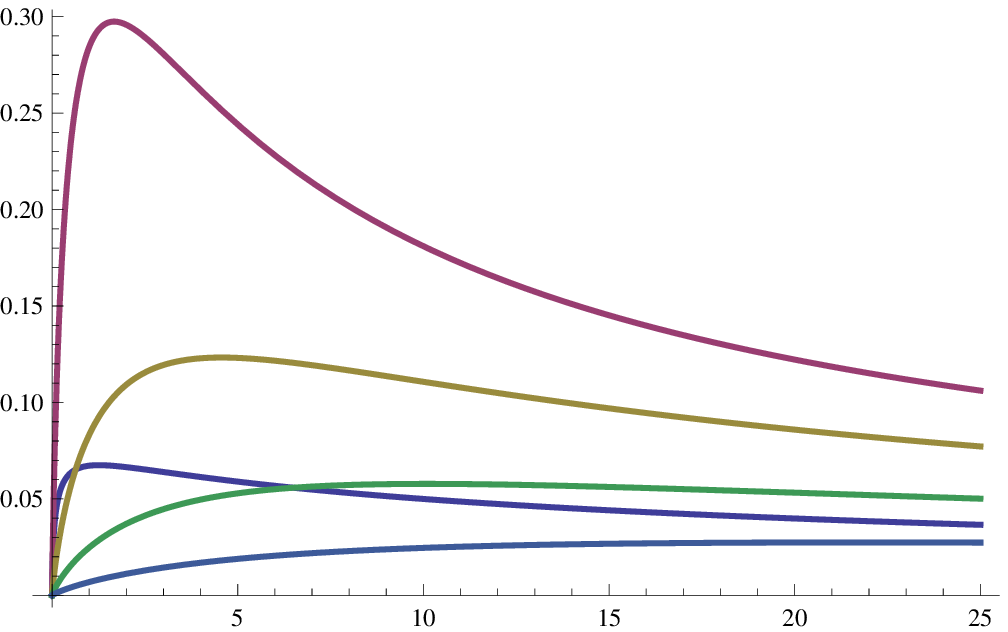}
}\end{center}}
\nnd{\bf Figure 5}\quad The difference ${\tilde k}_{q, p}-k_{q, p}$ for
$q=p+x\,(p=1.1,\, 2,\, 5,\, 10,\, 20)$ and $x$ varies over $(0.0001, 25)$. When $p=1.1$, the curve is special, located at
lower level and intersects with two others. The remaining curves from top to bottom
correspond to $p=2$, $5$, $10$, and $20$, respectively.
\medskip

\nnd Thus, our upper bound ${\tilde k}_{q,p} B$ given in (\ref{30}) is a little bigger than the basic one (\ref{22}).
As illustrated by Figures 1--4, $\dz_1$ improves $k_{q, p} B$
(not only ${\tilde k}_{q, p} B$) remarkably.
However, the proof for the sharp factor $k_{q, p}$\,(when $q> p$) becomes much more
technical (cf. \rf{cmp-07}{pages 45--47} for historical remarks and \ct{bga30}.
See also Example \ref{t2-5} below). Therefore, we prove only the upper bound given in (\ref{31}),
as ones often do so \rf{ok-90}{Theorem 1.14}.
Actually, one often regards (\ref{22}) replacing $k_{q, p}$ by ${\tilde k}_{q,p}$ as ``basic estimates'', due to the reasons just mentioned above.

Among $\dz_1$, ${\bar\dz}_1$, and ${\tilde \dz}_1$ in the corollary, the most complicated one is
${\tilde \dz}_1$. It is not simple even for the simplest Example \ref{t1-1}:
$$\aligned
&{\tilde \dz}_1=\sup_{z\in (0, 1]}\bigg(\frac{f_2^{(z)}(z)}{z}\bigg)^{(q-1)/p}
\big\|f_2^{(z)}\big\|_{\mu, q}^{1-q/ p},\\
& f_2^{(z)}(x)=\int_0^{x\wedge z} \d y \bigg[\frac 1 q (z^q-y^q)+z^{q-1}(1-z)\bigg]^{p^*-1},\qqd x\in [0, 1].
\endaligned$$
The main contribution of the sequence $\{{\tilde \dz}_n\}$ is, when $q=p$, its increasing property which then implies that $\{{\tilde \dz}_n\}$ is
closer and closer, step by step, to $A$. Therefore, the sequence $\{{\bar\dz}_n\}$ posses
the same property since ${\bar\dz}_{n+1}\ge {\tilde\dz}_n$ by \ct{cwz-13}.
However, there is no direct proof for the increasing property of the sequence
$\{{\bar \dz}_n\}$ even though it is believed to be true.
From \ct{cwz-13}, it is also known that
in the particular case of $q=p$, we have ${\bar\dz}_1\ge {\tilde \dz}_1$ if $p\ge 2$,
${\bar\dz}_1\le {\tilde \dz}_1$ if $p\in (1, 2]$, and ${\bar\dz}_1= {\tilde \dz}_1$ if $q=p=2$.
Thus, only in a small region of $(p, q)$, ${\tilde \dz}_1$ can be better than
${\bar\dz}_1$. For instance, setting $p=1.1$ in our Example \ref{t1-1}, then only for those $q\in [1.1, 1.55]$, one has ${\bar\dz}_1\le {\tilde \dz}_1$.
Next, let $p=2$, then we have ${\bar\dz}_1> {\tilde \dz}_1$ once $q>p$.
For this reason, unlike the case of $q=p$, here
we do not pay much attention to study the sequence $\{{\tilde \dz}_n\}$ in the case of $q>p$.

Having Corollary \ref{t2-3} at hand, it is not difficult to compute $\dz_1$ and ${\bar\dz}_1$
given in Example \ref{t1-1}. To obtain the constant $A^*$ there, we need more work.

\rmk\lb{t2-4}{\rm
We are now going to describe the upper estimate (\ref{19}) in a different way.
First, when $q=p$, we can rewrite $I\!I^*$ as $I\!I_{\mu, p}^{ \nu, p}$:
$$I\!I_{\mu, p}^{ \nu, p}(f)(x)=\frac{1}{f(x)}\int_0^x \d y \,{\hat v}(y) \bigg(\int_y^D f^{p-1}\d\mu\bigg)^{p^*-1}.$$
At the same time, we rewrite $A$ in (\ref{01}) as $A_{\mu, p}^{\nu, p}$.
In this case, in view of the first inequality of (\ref{19}), we have obtained
$$A_{\mu, p}^{\nu, p}\le \inf_{f\in {\scr F}_{I\!I}}\Big[\sup_{x\in (x, D)}I\!I_{\mu, p}^{ \nu, p}(f)(x)\Big]^{1/p^*}.$$
Actually, by \rf{cwz-13}{Theorem 2.1}, the equality sign here holds:
\be A_{\mu, p}^{\nu, p}= \inf_{f\in {\scr F}_{I\!I}}\Big[\sup_{x\in (x, D)}I\!I_{\mu, p}^{ \nu, p}(f)(x)\Big]^{1/p^*}.\lb{06}\de
Next, for general $p$ and $q$, we may use the similar notation
$I\!I_{\mu, q}^{\nu, p}$ and $A_{\mu, q}^{\nu, p}$. When $q\ge p$, noting that
corresponding to ${\tilde p}=q/p^*+1$ and ${\tilde v}=v^{q/p}$, we have
$$\hat{\tilde v}={\tilde v}^{-\frac{1}{\tilde p-1}}=\big(v^{\frac q p}\big)^{-\frac{1}{q/p^*+1-1}}
= v^{-\frac 1{p-1}}={\hat v},$$
It follows that
$I\!I^*=I\!I_{\mu,\, {\tilde p}}^{\nu^{q/p},\, {\tilde p}}$,
here $\nu^{q/p}$ denotes for a moment the measure when the density of $\nu$ is replaced by its power of $q/p$.
By using the first inequality of (\ref{19}) again, we have
$$\aligned
A_{\mu, q}^{\nu, p}&\le \inf_{f\in {\scr F}_{I\!I}}\Big[\sup_{x\in (x, D)} I\!I_{\mu, \, {\tilde p}}^{\nu^{q/p},\, {\tilde p}}(f)(x)\Big]^{1/p^*}\\
&=\bigg\{\inf_{f\in {\scr F}_{I\!I}}\Big[\sup_{x\in (x, D)}I\!I_{\mu,\, {\tilde p}}^{\nu^{q/p},\, {\tilde p}}(f)(x)\Big]^{{q}/{(q+p^*)}}\bigg\}^{{1}/{p^*}+ 1/q}
\endaligned$$
since the conjugate of ${\tilde p}=1+q/p^*$ is $1+p^*/q$. By (\ref{06}), we have thus obtain
the following estimate
\be A_{\mu, q}^{\nu, p}\le\Big[A_{\mu,\, {\tilde p}}^{\nu^{q/p}, \, {\tilde p}}\Big]^{{1}/{p^*}+ 1/ q}.\lb{35}
\de
In other words, when $q\ne p$, we are estimating the optimal constant $A_{\mu, q}^{\nu, p}$
of a mapping $L^p(\nu)\to L^q(\mu)$ by the one of
$L^{\tilde p}(\tilde\nu)\to L^{\tilde p}(\mu)$. When $q=p$, the right-hand side of
(\ref{35}) coincides with its left-hand side and so (\ref{35}) becomes an equality.
In Example \ref{t1-1}, the upper bound $A^*$ denotes the right-hand side of (\ref{35}).
Note that without assuming (\ref{06}), by part (1) of Theorem \ref{t2-1}, the estimate
(\ref{35}) is the best one we can expected. This indicates a limitation of (\ref{19}) since Figure 3 shows
that there is a small difference between the two sides of (\ref{35}) (see also Example \ref{t2-5} below). In contract to
part (1) of Theorem \ref{t2-1}, part (2) of the theorem can be sharp at least when there is a solution to the Euler-Lagrange equation (or ``eigenequation''):
$$(v {g'}^{p-1})' + u g^{q-1}=0, \qqd g, g'>0 \text{ on }(0, D).$$
}
\dermk

We conclude this section by looking an extremal example to which there is no room for improving
the upper estimate in (\ref{22}). Refer to Lemma \ref{ta-4} and Proposition \ref{ta-5} in Appendix for more
general results.

\xmp\lb{t2-5}{\rm Let $q>p>1$, $D=\infty$, $\mu(\d x)=x^{-q/p*-1}\d x$, and $\nu(\d x)=\d x$.
Then the optimal constant in the Hardy-type inequality is
$$A=\bigg(\frac{p^*}{q}\bigg)^{{1}/{q}}
\Bigg[\frac{\Gamma\Big(\frac{pq}{q-p}\Big)}{\Gamma\Big(\frac{q}{q-p}\Big) \Gamma\Big(\frac{p(q-1)}{q-p}\Big)}\Bigg]^{1/ p- 1/ q}
=\bigg(\frac{p^*}{q}\bigg)^{{1}/{q}} k_{q, p}$$
which can be attained by a simple optimizer $f$ having derivative
$$f'(x)=\frac{\az}{(\bz x^{\gz}+1)^{(\gz +1)/\gz}},\qqd \az, \bz>0,\; \gz=\frac{q}{p}-1.$$
Refer to \ct{bga30} or Appendix for details. Since
$$B=\sup_{x>0}{\hat\nu}(0, x)^{{1}/{p^*}}\mu(x, \infty)^{1/ q}
=\sup_{x>0} x^{{1}/{p^*}}\bigg[\int_x^{\infty}y^{-{q}/{p^*}-1}\bigg]^{1/ q}
=\bigg(\frac{p^*}{q}\bigg)^{{1}/{q}},$$
the upper bound
of the basic estimates in (\ref{22}) is sharp. Actually, this is where the optimal factor
$k_{q, p}$ comes from.

Even though there is now nothing more to do about the upper estimate of $A$, to understand what happened
in such an extremal situation, we compute $\dz_n$. Because $\fz(x)=x$, $f_1=\fz^{\gz^*}$,
where $\gz^*=\frac{q}{p^*+q}$, and
$$\aligned
&\int_y^D \fz^{q \gz^*/p^*}\d \mu=\int_y^{\infty} z^{q\gz^*/p^*-q/p^*-1}
=\frac{p^*}{q(1-\gz^*)} y^{q(\gz^*-1)/p^*},\\
&\int_0^x\!\! \d y \bigg[\!\int_y^D\!\! \fz^{q \gz^*/p^*}\d \mu\bigg]^{{p^*}/{q}}
\!\!\!=\!\bigg[\frac{p^*}{q(1-\gz^*)}\bigg]^{p^*/q}\!\!\!\int_0^x\!\! y^{\gz^*-1}\d y
\!=\! \frac{1}{\gz^*} \bigg[\frac{p^*}{q(1-\gz^*)}\bigg]^{p^*/q}\!x^{\gz^*}\!\!.\endaligned$$
we have
$$f_2(x)=\int_0^x \d y \bigg[\int_y^{\infty}f_1^{q/p^*}\d \mu\bigg]^{p^*/q}
=\frac{1}{\gz^*} \bigg[\frac{p^*}{q(1-\gz^*)}\bigg]^{p^*/q} f_1(x)=: C f_1(x).$$
By induction, it follows that $f_{n+1}=C^n f_1$ and hence
$$\dz_n=\bigg(\sup_x\frac{f_{n+1}(x)}{f_n(x)}\bigg)^{1/p^*}
=C^{1/p^*}=\bigg(1+\frac{p^*}{q}\bigg)^{1/p^*+1/q},\qqd n\ge 1.$$
It is now easy to check that $\dz_n={\tilde k}_{q, p}B\,(\ge {k}_{q, p}B)$ for all $n\ge 1$.
Thus, no improvement of the upper bound ${\tilde k}_{q, p}B$ can be made by our approach.
This is not surprising since $\dz_1$ is already a sharp estimate of the right-hand side
of (\ref{35}). To see this, let $q\downarrow p$, we get
$$A_{\mu, q}^{\nu, p}\to A_{\mu, p}^{\nu, p}=\bigg(\frac{p^*}{p}\bigg)^{1/p} k_{p, p}
=\bigg(\frac{p^*}{p}\bigg)^{1/p} {\tilde k}_{p, p}=p^*.$$
(Actually, when $q=p$, we come back to the original Hardy inequality, its optimal
constant is well known to be $p^*$.)
Then replacing $p$ with ${\tilde p}=q/p^*+1$, noting that ${\tilde p}^*=p^*/q+1$,
we obtain the optimal constant
$$A_{\mu,\, {\tilde p}}^{\nu^{q/p}\!,\, {\tilde p}}=p^*/q+1$$
appearing on the right-hand side of (\ref{35}) which is clearly equal to $\dz_n$.
For general $q> p$, $\dz_n\equiv {\tilde k}_{q, p}B$ is actually bigger than, and so can not improve ${k}_{q, p}B$.

Next, we compute ${\bar\dz}_1$. Because
$$\aligned
\frac{1}{\fz (x)^{q/p}}\int_0^x \fz^q\d\mu +\fz(x)^{q/p^*} \mu(x, D)
&= \frac{1}{x^{q/p}}\int_0^x y^{q/p-1} + x^{q/p^*} \int_x^{\infty} y^{-q/p^*-1}\\
&=\frac{p}{q}+\frac{p^*}{q} x^{q/p},
\endaligned$$
we have by Corollary \ref{t2-3},
$${\bar\dz}_1= \bigg[\frac{pp^*}{q}\bigg]^{1/q}$$
which is clearly bigger than $B$: ${\bar\dz}_1/B=p^{1/q}>1$, and hence improves the lower bound
of the basic estimates in (\ref{22}). Since ${\bar\dz}_1$ is not sharp, the lower bound can be
usually improved step by step using the sequence $\{{\bar\dz}_n\}$. By Corollary \ref{t2-3},
for general $q>p>1$, the ratio
${\dz_1}/{\bar\dz_1}$ is controlled by ${\tilde k}_{q, p}$. However, from our experience
we do have (without proof) that
$$\sup_{q>p>1}{\dz_1}/{\bar\dz_1}\le \sqrt{2}< 2=\sup_{q>p>1} k_{q, p}.$$
In this sense, the ratio of the estimates in (\ref{22}) is improved.

Some illustrations of $A\, (=k_{q,p}B)$ and its lower bound ${\bar\dz}_1$ are given in Figures 6 and 7.
From which, one sees that our estimates are still effective even in
such an extremal situation.

{\begin{center}{\includegraphics[width=11.0cm,height=7.5cm]{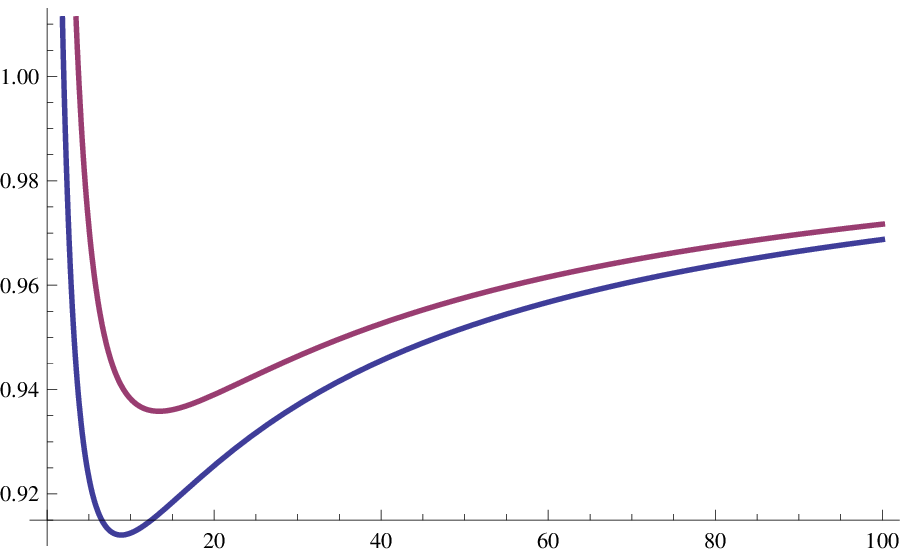}
}\end{center}}
\nnd{\bf Figure 6}\quad The constant $A=k_{q,p}B$ and its lower bound ${\bar\dz}_1$ in the case of $p=2$, $q=p+r$. $r\in [0, 100]$.

{\begin{center}{\includegraphics[width=11.0cm,height=7.5cm]{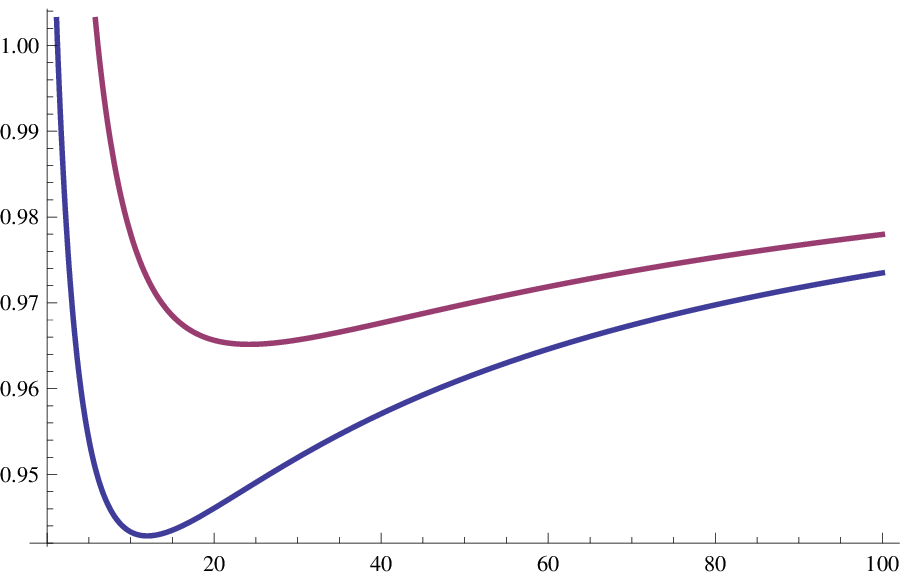}
}\end{center}}
\nnd{\bf Figure 7}\quad Everything is same as in the last figure except $p=2$ is replaced by $p=5$.
}
\dexmp

\section{Proofs}

It is now standard (cf. the explanation in the paragraph above \rf{cmf13}{(9)}
that to prove the main results stated in the last section, one
may assume that $\mu$ has a density $u$ if necessary. Similarly, one can
assume that $\nu$ has a density $v$. Besides, one can also assume some integrability
for ${\hat\nu}$ by an approximating procedure if necessary in the proofs below.
\medskip

\nnd{\bf Proof of Theorem \ref{t2-1}}. (a) First, we prove (\ref{19}).
Let $g$ satisfy $\|g\|_{\mu, q}=1$ and $g(0)=0$. Then for each positive $h$,
by a good use of the H\"older inequality, we have
\begin{align}
1&=\int_0^D g^q\d\mu=\int_0^D \mu(\d x)\bigg(\int_0^x g'\bigg)^q
= \int_0^D\mu(\d x) \bigg(\int_0^x g' v^{1/p} h^{-1} v^{-1/p} h\bigg)^q\nonumber\\
&\le \int_0^D\mu(\d x) \bigg(\int_0^x {g'}^p v h^{-p}\bigg)^{q/p}
     \bigg(\int_0^x {\hat v}  h^{p^*}\bigg)^{q/p^*}\qd(\text{since }p>1).\lb{36}
\end{align}
Here and in what follow, the Lebesgue measure $\d x$ is omitted.
To separate out the term $\int_0^D v {g'}^p$, we need an exchange of the order
of integration. When $q=p$, this is not a problem: one
simply uses Fubini's theorem and nothing is lost. However, when
$q\ne p$, this is not trivial. Fortunately, for $q\ge p$, one
can apply the H\"older-Minkowski inequality:
$$\aligned
&\bigg\{\int_{E_1} \mu (\d x)
\bigg[ \int_{E_2} f(x,y) \nu (\d y)\bigg]^{r} \bigg\}^{1/r}
\le \int_{E_2} \nu (\d y) \bigg[ \int_{E_1} f(x,y)^r  \mu (\d x)\bigg]^{1/r},\\
&\qd \mu, \nu: \sz\text{-finite measures, } r\in [1, \infty),\; f\ge 0.
\endaligned$$
Applying this inequality to $r=q/p$, $E_1=E_2=[0, D]$,
$\nu (\d y)=\big({g'}^p v h^{-p}\big)(y)\d y$, and
$$f(x, y)=\mathbbold{1}_{[0, x]} (y)\bigg(\int_0^x {\hat v}  h^{p^*}\bigg)^{p/p^*},$$
it follows that the right-hand of (\ref{36}) is controlled by
$$\bigg\{\int_{0}^D \d y \big({g'}^p v h^{-p}\big)(y)
\bigg[\int_y^D \mu (\d x) \bigg(\int_0^x {\hat v}  h^{p^*}\bigg)^{q/p^*}\bigg]^{p/q}\bigg\}^{q/p}.$$
Note that here we have only ``$\le $'' rather than ``$=$''.
Now, making a power $1/q$, we get
$$\aligned
1&\le \bigg\{\int_{0}^D \d y  \big({g'}^p v h^{-p}\big)(y)
\bigg[\int_y^D \mu (\d x) \bigg(\int_0^x {\hat v}  h^{p^*}\bigg)^{q/p^*}\bigg]^{p/q}\bigg\}^{1/p}\\
&\le \bigg(\int_{0}^D {g'}^p v\bigg)^{1/p}
\bigg\{\sup_{y\in (0, D)} \frac{1}{h(y)}
\bigg[\int_y^D \mu(\d x) \bigg(\int_0^x {\hat v} h^{p^*}\bigg)^{q/p^*}\bigg]^{1/q}\bigg\}.
\endaligned$$
Replacing $h$ by $h^{1/q}$, it follows that
\begin{equation} 1\le \bigg(\int_{0}^D {g'}^p v\bigg)^{1/p}
\bigg\{\sup_{y\in (0, D)} \frac{1}{h(y)}
\int_y^D \mu(\d x) \bigg(\int_0^x {\hat v} h^{p^*/q}\bigg)^{q/p^*}\bigg\}^{1/q}.\lb{37}
\end{equation}

To move further, we need an extension of the mean value theorem for integrals.

\lmm\lb{t3-1}{\cms Let $g>0$ on $(\az, \bz)$ and $\int_{\az}^{\bz}g\d\mu<\infty$.
Suppose that the integral $\int_{\az}^{\bz}f \d\mu$ exists (may be $+\infty$). Then
$$\sup_{x\in (\az, \bz)} \frac{\int_x^{\bz}f\d\mu}{\int_x^{\bz}g\d\mu}\le \sup_{x\in (\az, \bz)} \frac{f}{g}(x)
\qd\text{\cms and dually}\qd
\inf_{x\in (\az, \bz)} \frac{\int_x^{\bz}f\d\mu}{\int_x^{\bz}g\d\mu}\ge \inf_{x\in (\az, \bz)} \frac{f}{g}(x).
$$
}
\delmm

\prf Set
$\xi= \sup_{x\in (\az, \bz)} \big( f/ g\big)(x).$
Without loss of generality, assume that $\xi<\infty$. Otherwise, the first
assertion is trivial. By assumptions, $g>0$ and moreover
$$f\le \xi g\qd\text{on}\qd {(\az, \bz)}.$$
Making integration over the interval $(x, \bz)$, it follows that
$$\int_{x}^{\bz} f\d\mu\le \xi \int_{x}^{\bz} g\d\mu, \qqd x\in (\az, \bz).$$
The first assertion then follows since $\int_{\az}^{\bz}g\d\mu\in (0, \infty)$.
Dually, we can prove the second assertion.
\deprf

We now come back to the proof of the inequality in (\ref{19}). Actually, we prove a
(formally) stronger conclusion. Let
$${\scr F}_{I\!I}^*=\{f: f(0)\ge 0,\; f>0 \text{ on }(0, D)\}\qd\big[\supset {\scr F}_{I\!I}\big].$$
For a given $f\!\in\! {\scr F}_{I\!I}^*$, without loss of generality, assume that
$\sup_{x\in (0, D)}I\!I^*(f)(x)$ $<\infty$. Otherwise, the upper bound we are going to
 prove is trivial.
Let $h(x)=\int_x^D f^{q/p^*}\d\mu$. As an application
of Lemma \ref{t3-1}, since $ h<\infty$, we have
$$\aligned
&\sup_{x\in (0, D)} \frac{1}{h(x)}
\int_x^D \mu(\d y) \bigg(\int_0^y {\hat v} h^{p^*/q}\bigg)^{q/p^*}\\
&\qd\le \bigg\{\sup_{x\in (0, D)}\frac{1}{f(x)}\int_0^x \d y \, {\hat v}(y) \bigg[\int_y^D f^{q/p^*}\d\mu\bigg]^{p^*/q} \bigg\}^{q/p^*}.
\endaligned$$
Inserting this into (\ref{37}) and making supremum with respect to $g$, it follows that
$$ A\le \Big[\sup_{x\in (0, D)}I\!I^*(f)(x)\Big]^{1/p^*}$$
and then
$$ A \le \inf_{f\in {\scr F}_{I\!I}^*}\Big[\sup_{x\in (0, D)}I\!I^*(f)(x)\Big]^{1/p^*}
\le \inf_{f\in {\scr F}_{I\!I}}\Big[\sup_{x\in (0, D)}I\!I^*(f)(x)\Big]^{1/p^*}.$$
This gives us the first inequality in (\ref{19}).
Furthermore, applying Lemma \ref{t3-1} again, we obtain
$$A\le \inf_{f\in {\scr F}_{I\!I}}\Big[\sup_{x\in (0, D)}I\!I^*(f)(x)\Big]^{1/p^*}
\le \inf_{f\in {\scr F}_{I}}\Big[\sup_{x\in (0, D)}I^*(f)(x)\Big]^{1/p^*}.$$
Now, for a given $f\in {\scr F}_{I\!I}$ with $\sup_x I\!I^*(f)<\infty$, let
$g=fI\!I^*(f)$. Then $g\in {\scr F}_I$ and
$$g'(x)={\hat v}(x)\bigg[\int_x^D f^{q/p^*}\d \mu\bigg]^{p^*/q}
\ge  {\hat v}(x)\bigg[\int_x^D g^{q/p^*}\d \mu\bigg]^{p^*/q}\Big[\inf_x I\!I^*(f)^{-1}\Big].$$
That is,
$$\sup_x I\!I^*(f)(x)\ge \frac{{\hat v}(x)}{g'(x)}\bigg[\int_x^D g^{q/p^*}\d \mu\bigg]^{p^*/q},$$
and then
$$\sup_x I\!I^*(f)(x)\ge \sup_x I^*(g)(x).$$
On both sides, making successively, power $1/p^*$, infimum with respect to $g\in {\scr F}_I$,
and then infimum with respect to $f\in {\scr F}_{I\!I}$, we obtain
$$\inf_{f\in {\scr F}_{I\!I}}
\bigg[\sup_x I\!I^*(f)(x)\bigg]^{1/p^*}
\ge \inf_{g\in {\scr F}_{I}} \bigg[\sup_x I^*(g)(x)\bigg]^{1/p^*}.$$
Therefore, the equality in (\ref{19}) holds.

(b) Next, we prove (\ref{20}). Given $f\in {\widetilde {\scr F}}_{I\!I}$, define $g_0=[f I\!I(f)](\cdot\wedge x_0)$. Then
$$\int_0^D v {g_0'}^p=\int_0^D v {g_0'}^{p-1}\d g_0
=\big(v g_0{g_0'}^{p-1}\big)(x_0-)-\int_0^{x_0} g_0 \big(v {g_0'}^{p-1}\big)'.$$
By definition of $g_0$, we have
$$\big(v g_0{g_0'}^{p-1}\big)(x_0-)=g_0(x_0)\int_{x_0}^D f^{q-1}\d\mu,\qqd
\big(v {g_0'}^{p-1}\big)'=- f^{q-1}\d\mu.$$
Hence we have
$$\int_0^D v {g_0'}^p=\int_0^D g_0 f^{q-1}\d\mu
\le \bigg(\sup_{x\in (0, D)}\frac f {g_0}\bigg)^{q-1} \int_0^D g_0^q\d\mu.$$
That is,
$$\|g'\|_{\nu, p}\le \bigg(\sup_{x\in (0, D)}\frac f {g_0}\bigg)^{(q-1)/p}\|g_0\|_{\mu, q}^{q/p}.$$
In other words,
$$\frac{\|g_0'\|_{v, p}}{\|g_0\|_{\mu, q}}\le \bigg(\sup_x \frac f {g_0}\bigg)^{(q-1)/p} \|g_0\|_{\mu, q}^{q/p-1}.$$
We have thus obtain
$$A\ge \sup_{f\in {\widetilde{\scr F}}_{I\!I}} \|f I\!I\, \tilde{}(f)\|_{\mu, q}^{1-q/p}\Big(\inf_{x\in (0, D)} I\!I\, \tilde{}(f)(x)\Big)^{(q-1)/p}.$$
This proves the first assertion of part (2). Then the second one follows by using
the proof similar to the last part of proof (a).
\deprf

\medskip

\nnd{\bf Proof of Theorem \ref{t2-2}}. The approximating sequences $\{\dz_n\}$ and $\{{\tilde\dz}_n\}$ are simply
successive application of Theorem \ref{t2-1}. The sequence $\{{\bar\dz}_n\}$ is a direct application of (\ref{08}).
The monotonicity of $\dz_n$ in $n$ is obtained by using Lemma \ref{t3-1} twice.
\deprf

To prove the (basic) upper bound given in Corollary \ref{t2-3}, we need the following result.

\lmm\lb{t3-2}{\cms Let $\fz>0 $ on $(0, D)$ and
$$B:=\sup_{x\in (0, D)} \fz (x)^{1/p^*} \mu(x, D)^{1/q}<\infty.$$
Then for each $\gz\in (0, 1)$, we have
$$\bigg(\int_x^D \fz^{\gz q/p^*}\d\mu\bigg)^{1/q}\le
\frac{B}{(1-\gz)^{1/q}}\fz^{(\gz-1)/p^*}.$$
}
\delmm

\prf For a function $h\in C[0, D]\cap C^1(0, D)$ with $h(0)=0$, write
$$\int_x^D h\d\mu=-\int_x^D h(y) \d M(y),\qqd M(y):=\mu(y, D).$$
Applying \rf{cmf00}{Proof of Lemma 1.2} to
$c=B^q$ with a change of $\fz$ by $\fz^{q/p^*}$, it follows that
$$\int_x^D \fz^{\gz q/p^*}\d\mu\le \frac{B^q}{1-\gz} \fz^{q(\gz-1)/p^*}. $$
The required assertion now follows immediately.
\deprf

\nnd{\bf Proof of Corollary \ref{t2-3}}. The main assertion as well as the
formula of $\dz_1$
are obtained by Theorem \ref{t2-2} directly, except the estimates involving $B$
and the formulas of ${\bar\dz}_1$ and ${\tilde\dz}_1$. The inequality involving
$B$ in the middle is based on (\ref{22}).

(a) To prove the upper bound given in (\ref{31}), we specify $\fz$ used in Lemma \ref{t3-2}: $\fz(x)={\hat \nu}(0, x)$, and set $f=\fz^{\gz}$. Then
$$\frac{\hat v}{f'}=\frac{1}{\gz \fz^{\gz -1}}.$$
By Lemma \ref{t3-2}, we have
$$\big[I^*(f)(x)\big]^{1/p^*}\le
B \gz^{-1/p^*} (1-\gz)^{-1/q}.$$
Optimizing the right-hand side with respect to $\gz$, the minimum
$$ \bigg(1+\frac{q}{p^*}\bigg)^{1/q}
 \bigg(1+\frac{p^*}{q}\bigg)^{1/p^*}\qd \big[={\tilde k}_{q,p}\big]$$
of $\gz^{-1/p^*} (1-\gz)^{-1/q}$ is attained at
$$\gz^*= \frac{q}{p^*+q}.$$
We have finally arrived at $\dz_1\le {\tilde k}_{q, p} B$ by using the equality in
(\ref{19}) with the specific $f=\fz^{\gz^*}$. From the proof, the main reason why
$\dz_1$ can improve ${\tilde k}_{q, p} B$ is clear: $\dz_1$ is defined by using the
operator $I\!I^*$, but its upper bound ${\tilde k}_{q, p} B$ is deduced from the
operator $I^*$. Usually, there is a gap between $\sup_x I\!I^*(f)$ and $\sup_x I^*(f)$
for a fixed $f$.

(b) To compute ${\bar\dz}_1$, recall our test function
$$\fz(x)=\fz^{(x_0)}(x)=\int_0^{x\wedge x_0} {\hat v}.$$
The reason to choose this function is the following observation:
$$\int_0^D v {\fz'}^p=\int_0^{x_0}v^{p(1-p^*)} v=\int_0^{x_0}{\hat v}=\fz(x_0).$$
Next, because of
$$\int_0^D \fz^q\d\mu=\int_0^{x_0} \fz^q\d\mu+\fz(x_0)^q \mu(x_0, D),$$
it follows that
$$\frac{\|\fz\|_{\mu, q}}{\|\fz'\|_{v, p}}=
\bigg[\frac{1}{\fz(x_0)^{q/p}}\int_0^{x_0} \fz^q\d\mu+\fz(x_0)^{q/p^*} \mu(x_0, D)\bigg]^{1/q}.$$
Making supremum with respect to $x_0$, we obtain ${\bar\dz}_1$.

The proof of ${\bar\dz}_1\ge B$ is rather easy. Simply ignore the first term
in the sum in (\ref{32}). The improvement of ${\bar\dz}_1$ from $B$ is obvious.

(c) To compute ${\tilde \dz}_1$, recall that
$$\aligned
f_1^{(x_0)}(x)&={\hat\nu} (0, x\wedge x_0),\\
f_2^{(x_0)}(x)\!&=\!
\!\int_0^{x\wedge x_0}\! \d y\,{\hat v}(y)\bigg[\int_y^{x_0}\!\! \fz^{q-1}\d\mu
   \!+\! \fz(x_0)^{q-1}\mu (x_0, D)\bigg]^{p^*-1}\!\!,\qd x\in [0, D].
\endaligned
$$
For simplicity, in what follows, we ignore the superscript $(x_0)$ in $f_1^{(x_0)}$
and $f_2^{(x_0)}$. Clearly, we have
$$\inf_{x\in (0, D)} f_2(x)/f_1(x)= \inf_{x\in (0, x_0)}f_2(x)/f_1(x)$$
by the convention that $1/0=\infty$.
Next, we show that the derivative of $f_2/f_1$ is non-positive
on $(0, x_0)$, that is
$$\fz(x)\bigg[\int_x^D \fz(\cdot \wedge x_0)^{q-1}\d\mu\bigg]^{p^*-1}
-\int_0^x \d y\,{\hat v}(y)\bigg[\int_y^D \fz(\cdot \wedge x_0)^{q-1}\d\mu\bigg]^{p^*-1}\le 0$$
on $(0, x_0)$.
This is obvious since for each $h\ge 0$,
$$\int_0^x \d y\,{\hat v}(y)\bigg[\int_y^D h\d\mu\bigg]^{p^*-1}
\ge \bigg(\int_0^x {\hat v}\bigg) \bigg[\int_x^D h\d\mu\bigg]^{p^*-1}
=\fz(x)\bigg[\int_x^D h\d\mu\bigg]^{p^*-1}.$$
Hence we indeed have
$$\inf_{x\in (0, D)} f_2(x)/f_1(x)= f_2(x_0)/f_1(x_0).$$
We have thus obtained ${\tilde\dz}_1$ as stated in the corollary.
\deprf

\section{Appendix. The inequalities on finite intervals and the sharp factor}

As far as we know, the basic estimates (\ref{02}) with universal optimal constant $k_{q, p}$
was proved only for the half-line (cf. \ct{mvm92}). In this appendix, we show that the
estimates with the same factor $k_{q, p}$ actually hold for every finite interval.
The study on this problem also provides us a chance to examine how to obtain (\ref{02}).
The main result of this section is Theorem \ref{ta-7}.
We begin with our study on three comparison results for the optimal constants
and their basic upper estimates in different intervals. The first one is a comparison for the optimal constants only.

\lmm\lb{ta-1}{\cms Let $A_D$ be the optimal constant in the Hardy-type inequality on the interval
$(0, D)$. Then we have $A_{D}\uparrow A_{D'}$ as $D\uparrow D'\le \infty$.
Here we use the same notation $(\mu, \nu)$ to denote the Borel measures on
$[0, D']$ and their restriction to $[0, D]$. In particular, if the inequality holds on $(0, D')$,
then it also holds with the same constant $A_{D'}$ on $(0, D)$ for every $D<D'$.}
\delmm

\prf (a) Extending $f$ from $[0, D]$ to $[0, D')$ by setting
$f=f(\cdot \wedge D)$, it follows that
$$\bigg[\!\int_0^D\! |f|^q\d \mu\bigg]^{1/q}
\!\!\le\! \!\bigg[\!\int_0^{D'}\! |f|^q\d \mu\bigg]^{1/q}
\!\!\le\! A_{D'}\!\bigg[\!\int_0^{D'}\! |f'|^p\d \nu\bigg]^{1/p}
\!\!=\!A_{D'}\!\bigg[\!\int_0^D\! |f'|^p\d \nu\bigg]^{1/p}\!.$$
The last assertion of the lemma is now obvious. We have thus proved the monotonicity:
$A_D\le A_{D'}$ whenever $D\le D'$.

(b) To prove the convergence in the first assertion, consider first the simplest case that
$\mu[0, D']=\infty$. Then $D'=\infty$ since $\mu$ is Borel. Clearly, we have $B_{D'}=\infty$
and so is $A_{D'}$ by our basic estimates. Besides, restricting to $[0, n]$, we have
$$A_n\ge B_n=\sup_{x\in [0, n]}\mu[x, n]^{1/q}{\hat\nu}[0, x]^{1/p^*}
\ge \mu[1, n]^{1/q}{\hat\nu}[0, 1]^{1/p^*} \to \infty\;\text{as }n\to\infty,$$
hence the convergence in the first assertion holds in this case.

(c) Let $\mu[0, D']<\infty$ and $A_{D'}<\infty$.
Then for every $f$ satisfying $\|f'\|_{\nu, p}\in (0, \infty)$, we have
$$\frac{\|f \mathbbold{1}_{[0, D]}\|_{\mu, q}}{\|f' \mathbbold{1}_{[0, D]}\|_{\nu, p}}\to
\frac{\|f\|_{\mu, q}}{\|f'\|_{\nu, p}}\le A_{D'}\qqd\text{as } D\uparrow D'.$$
 Since $A_{D'}<\infty$, for every $\vz>0$, we can choose first $f=f_{\vz}$ such that $\|f'\|_{\nu, p}\in (0, \infty)$ and
$$A_{D'}\le \frac{\|f\|_{\mu, q}}{\|f'\|_{\nu, p}}+\vz,$$
then we can choose $D$ closed to $D'$ such that
$$\frac{\|f\|_{\mu, q}}{\|f'\|_{\nu, p}}\le\frac{\|f \mathbbold{1}_{[0, D]}\|_{\mu, q}}{\|f' \mathbbold{1}_{[0, D]}\|_{\nu, p}}  +\vz.$$
Therefore, we obtain
$$A_D\le A_{D'}\le \frac{\|f \mathbbold{1}_{[0, D]}\|_{\mu, q}}{\|f' \mathbbold{1}_{[0, D]}\|_{\nu, p}}  +2\vz\le A_D+2\vz.$$
From this, we conclude that the convergence also holds in the present case.

(d) Finally, the proof in the case that $\mu[0, D']<\infty$ but $A_{D'}=\infty$ is in parallel to the proof (c).
\deprf

The next result is a comparison of the factor in the basic estimates for
different intervals.

\lmm\lb{ta-2}{\cms Let $A_D(\mu, \nu)$ and $B_D(\mu, \nu)$ denote the constants $A$ and $B$, respectively, given in the basic estimates (\ref{02}) for the inequality on interval $[0, D]$ with measures $\mu$ and $\nu$.
Next, let $D<D'\le \infty$ and $(\mu', \nu')$ be an extension of $(\mu, \nu)$ to $[0, D')$: $\mu'|_{[0, D]}=\mu$, $\nu'|_{[0, D]}=\nu$, and moreover $\mu'|_{(0, D')}=0$.
\begin{itemize} \setlength{\itemsep}{-0.8ex}
\item[(1)]
Suppose that $A_{D'}(\mu', \nu')\le k B_{D'}(\mu', \nu')$ for a universal constant $k$, then we have $A_D(\mu, \nu)\le k B_{D}(\mu, \nu)$.
\item[(2)] In particular, if the inequality in part (1) holds for arbitrary (resp. absolutely continuous) pair $(\mu', \nu')$, then so does the conclusion for
arbitrary (resp. absolutely continuous) pair $(\mu, \nu)$.
\end{itemize}}
\delmm

\prf Clearly, we need only to prove the first assertion. Then the second one
follows immediately.
As in the last proof, extend $f$ from $[0, D]$ to $[0, D')$ by setting $f=f(\cdot \wedge D)$.
Then we have
$$\aligned
\bigg[\int_0^D |f|^q\d \mu\bigg]^{1/q}
&=\bigg[\int_0^{D'} |f|^q\d \mu'\bigg]^{1/q}\qd\text{\big(since $\mu'|_{(D, D')}=0$\big)}\\
&\le A_{D'}(\mu', \nu') \bigg[\int_0^{D'} |f'|^p\d \nu'\bigg]^{1/p}\qd\text{(by definition of $A_{D'}(\mu', \nu')$)}\\
&\le k B_{D'}(\mu', \nu') \bigg[\int_0^{D'} |f'|^p\d \nu'\bigg]^{1/p}\qd\text{(by assumption)}\\
&= k B_{D'}(\mu', \nu') \bigg[\int_0^{D} |f'|^p\d \nu\bigg]^{1/p}\qd\text{\big(since $f'|_{(D, D')}=0$\big)}.
\endaligned$$
Because
$$\aligned
B_{D'}(\mu', \nu')&=\sup_{x\in (0, D')} {\hat\nu'}(0, x)^{1/p^*}\mu'(x, D')^{1/q}\\
&=\sup_{x\in (0, D)} {\hat\nu'}(0, x)^{1/p^*}\mu'(x, D)^{1/q}\qd\text{\big(since $\mu'|_{(D, D')}=0$\big)}\\
&=\sup_{x\in (0, D)} {\hat\nu}(0, x)^{1/p^*}\mu(x, D)^{1/q}\qd\text{\big(since $\mu'|_{[0, D)}=\mu$ and $\nu'|_{[0, D]}=\nu$\big)}\\
&=B_D (\mu, \nu),
\endaligned$$
it follows that
$$\bigg[\int_0^D |f|^q\d \mu\bigg]^{1/q}
\le k B_D (\mu, \nu) \bigg[\int_0^{D} |f'|^p\d \nu\bigg]^{1/p}.$$
Hence $A_D(\mu, \nu)\le k B_D (\mu, \nu)$ as required.
\deprf

The next result is somehow a refinement of Lemma \ref{ta-1}, but in an opposite way: from local sub-intervals to the whole interval. It provides us
an approximating procedure for unbounded interval.

\lmm\lb{ta-6}{\cms Given Borel measures $\mu^D$ and $\nu^D$ on $[0, D]$, extend them to $[0, {D'})$, $D<{D'}\le\infty$, as follows:
$${\tilde\mu}^D=\begin{cases}
\mu^D\qd  &\text{\cms on } [0, D],\\
0 \qd &\text{\cms on } (D, {D'});
\end{cases}$$
$${\tilde\nu}^{D, \#}=\begin{cases}
\nu^{D}\qd  &\text{\cms on } [0, D],\\
\# \qd &\text{\cms on } (D, {D'}),
\end{cases}$$
where $\#$ is an arbitrary Borel measure. Then we have
$A_D=A\big({\tilde\mu}^D, {\tilde\nu}^{D, \#}\big)$ and
$B_D=B\big({\tilde\mu}^D, {\tilde\nu}^{D, \#}\big)$.
}\delmm

\prf Following the proof of Lemma \ref{ta-2}, it is easy to check that
$B_D=B\big({\tilde\mu}^D, {\tilde\nu}^{D, \#}\big)$. Next,
applying the inequality
$$\|f\|_{L^q({\tilde\mu}^D)}\le A\big({\tilde\mu}^D, {\tilde\nu}^{D, \#}\big)\big\|f'\big\|_{L^p({\nu}^{D, \#})},$$
to $f^D=f(\cdot \wedge D)$, we obtain
$$\big\|f^D\big\|_{L^q({\mu}^D)}\le A\big({\tilde\mu}^D, {\tilde\nu}^{D, \#}\big)\big\|\big(f^D\big)'\big\|_{L^p({\nu}^{D})}.$$
Because $f$ is arbitrary and so is $f^D$, this implies that $A_D\le A\big({\tilde\mu}^D, {\tilde\nu}^{D, \#}\big)$. Conversely,
for every function $f$ on $(0, D')$ with $f(0)=0$, we have
$$\aligned
\bigg[\int_0^{D'} |f|^q\d {\tilde\mu}^D \bigg]^{1/q}
&=\bigg[\int_0^D |f|^q\d {\mu}^D \bigg]^{1/q}\\
&\le A_D \bigg[\int_0^D |f'|^p\d {\nu}^D \bigg]^{1/p}\;\text{(by definition of $A_D$)}\\
&\le A_D \bigg[\int_0^{D'} |f'|^p\d {\tilde\nu}^{D, \#} \bigg]^{1/p}\;\text{(since $D<D'$)}.
\endaligned$$
This implies that $A\big({\tilde\mu}^D, {\tilde\nu}^{D, \#}\big)\le A_D$ and
then the equality holds.
\deprf

\lmm[Bliss, 1930]\lb{ta-4}{\cms Let $q>p\,(p, q\in (1, \infty))$, $\nu(\d x)=\d x$, and $\mu (\d x)=x^{-q/p^*-1}\d x$ on $[0, D]$. Then we have
 $A\le k_{q, p}\, (p^*/q)^{1/q}$ with equality sign holds provided $D=\infty$.}
\delmm

\prf The case that $D=\infty$ was proved in Bliss' original paper \ct{bga30}. Then by Lemma \ref{ta-1}, the conclusion also holds for finite $D$.\deprf

The next result is a generalization of Bliss's lemma. It says that the basic
upper estimate in (\ref{02}) is sharp for a large class of $(\mu, \nu)$.

\prp\lb{ta-5}{\cms Let $q>p\,(p, q\in (1, \infty))$, $\nu(\d x)=v(x)\d x$, and define
${\hat v}(x)=v(x)^{1/(1-p)}$. Then the Hardy-type inequality holds
on $[0, D]$ with $\mu(\d x):=u(x)\d x$,
$$0\le u(x) \le -B_1^q \; \frac{\d}{\d x}\bigg(\int_0^x \hat v\bigg)^{-q/p*},$$
where $B_1\in (0, \infty)$ is a constant. Moreover, its optimal constant $A_D$ satisfies $A_D\le k_{q, p} B_1$.
In particular, when $D=\infty$,
\be{\hat\nu}(0, \infty)=\infty\quad \text{\cms and} \quad
\sup_{x\in (0, D)}\bigg[\int_0^x \hat v\bigg]^{1/p*}\bigg[\int_x^\infty u\bigg]^{1/q}
= B_1,  \lb{a-07}\de
the upper bound is sharp with $B_1=B$ defined by (\ref{23}).}\deprp

\prf Throughout the proof, we restrict ourselves to the special case that
$$u(x) = -B_1^q \; \frac{\d}{\d x}\bigg(\int_0^x \hat v\bigg)^{-q/p^*}>0.$$
The general case stated in the proposition then follows immediately. In this situation,
the last assertion of the proposition is due to \rf{mvm92}{Theorem 1}. Actually, the essential part of the proof is in the special case that $B_1=1$. The use of $B_1$ indicates an additional freedom for the choice of $u$, even in the present non-linear situation.

(a) By definition of $u$, we have
$$\int_x^D u= B_1^q\bigg[\bigg(\int_0^x \hat v\bigg)^{-q/p*}
- \bigg(\int_0^D \hat v\bigg)^{-q/p*} \bigg]
\le  B_1^q \bigg(\int_0^x \hat v\bigg)^{-q/p*}.$$
Note that here the equality sign holds iff so does the first condition in (\ref{a-07}).
Hence
$$B=\sup_{x\in (0, D)}\bigg[\int_0^x \hat v\bigg]^{1/p*}\bigg[\int_x^D u\bigg]^{1/q}
\le B_1<\infty,$$
and $B_1=B$ once the first condition in (\ref{a-07}) holds. Then the second
condition in (\ref{a-07}) is automatic in the present special case.

(b) Define
\bg{align}
s(x)&=\int_0^x {\hat v},  \lb{a-08}\\
\fz(s(x))&= f(x){\hat v}(x)^{-1}\qd (\fz=\fz_f^{}).  \lb{a-09}
\end{align}
Since $s(x)$ is increasing in $x$, its inverse function $s^{-1}$ is well-defined.
Then the last equation can be rewritten as
$$\fz (s)= f\big(s^{-1}(s)\big){\hat v}\big(s^{-1}(s)\big)^{-1}.$$
Because
\be f(x)\d x= f(x) {\hat v}(x)^{-1}\d s(x)=\fz (s)\d s,\lb{a-10}\de
we have by (\ref{a-10}),
$$H f(x):=\int_0^x f=\int_0^{s(x)} \fz=H\fz (s(x)).$$
Next, because of definition of $u$ and $s$,
\bg{gather}
u(x)= \frac{q}{p^*} B_1^q \bigg(\int_0^x {\hat v}\bigg)^{-q/p^*-1} {\hat v}(x)
= \frac{q}{p^*} B_1^q s(x)^{-q/p^*-1} {\hat v}(x),\nonumber\\
u(x)\d x = \frac{q}{p^*} B_1^q s(x)^{-q/p^*-1} \d s(x)\qd\text{(cf. (\ref{a-08}))},\nonumber
\end{gather}
we obtain
$$\aligned
\bigg[\int_0^D \big(Hf(x)\big)^q u(x)\d x\bigg]^{1/q}
&=B_1\bigg(\frac{q}{p^*}\bigg)^{1/q}
\bigg[\int_0^{s(D)} \big(H\fz(s)\big)^q\, s^{-q/p^*-1}\d s\bigg]^{1/q}\\
&\le k_{q, p} B_1\bigg[\int_0^{s(D)} \fz(s)^p\d s\bigg]^{1/p}
\endaligned$$
by Bliss' lemma and Lemma \ref{ta-1}. The equality sign holds once $s(D)=\infty$, i.e.
(\ref{a-07}) holds. Since by (\ref{a-09}),
$$f(x)^p v(x)=\fz\big(s(x)\big)^p {\hat v}(x)^p v(x)=\fz\big(s(x)\big)^p {\hat v}(x),$$
and then by (\ref{a-08}),
$$f(x)^p v(x)\d x=\fz\big(s(x)\big)^p\d s(x), $$
we have
$$\bigg[\int_0^{s(D)} \fz(s)^p\d s\bigg]^{1/p}=\bigg[\int_0^D f(x)^p v(x)\d x\bigg]^{1/p}.$$
Therefore, we have proved that
$$\bigg[\int_0^D \big(Hf(x)\big)^q u(x)\d x\bigg]^{1/q}
\le k_{q, p} B_1\bigg[\int_0^D f(x)^p v(x)\d x\bigg]^{1/p}.$$
This leads to the conclusion that $A_D\le k_{q, p} B_1$ as required. Again,
the equality sign holds under (\ref{a-07}).
\deprf

We can now state the main result in this section. When $D=\infty$, it is just \rf{bg91}{Theorem 8}.
If additionally (\ref{a-07}) holds, then it is \rf{mvm92}{Theorem 2}.

\thm\lb{ta-7}{\cms Let $q>p\,(p, q\in (1, \infty))$ and $D\le\infty$.
Then the basic estimates in (\ref{02}) hold for given $\mu$ and $\nu$.}\dethm

\prf The lower estimate in (\ref{02}) is shown in the proof of (\ref{30})(Corollary \ref{t2-3}). Our main task
is to prove the upper estimate in (\ref{02}).

Without loss of generality, assume that $\mu(\d x)=u(x)\d x$ and $\nu(\d x)=v(x)\d x$ on $[0, D]$
(cf. \ct{mb72}), and moreover $B<\infty$. Next, by part (2) of Lemma \ref{ta-2},
it suffices to prove the case that $D=\infty$. Note that
\begin{align}
\int_0^\infty \big(H f(x)\big)^q u(x)\d x
&=\int_0^\infty \bigg(\int_0^x\d \big(H f(t)\big)^q \bigg) u(x)\d x\nonumber\\
&=\int_0^\infty \bigg(\int_t^\infty u(x)\d x\bigg) \d \big(H f(t)\big)^q\qd\text{(by Fubini's theorem)}\nonumber\\
&\le B^q \int_0^\infty \bigg(\int_0^t {\hat v}(x)\d x\bigg)^{-q/p^*} \d \big(H f(t)\big)^q\nonumber\\
&\qqd\text{(by definition of $B$)}\nonumber\\
&= B^q \int_0^\infty s(t)^{-q/p^*} \d \big(H f(t)\big)^q\qd\text{(by (\ref{a-08}))}. \lb{a-11}
 \end{align}
 Next, note that
 $$\aligned
 \int_t^\infty s(x)^{-q/p^*-1}\d s(x)
 &= -\frac{p^*}{q}s(x)^{-q/p^*}\big|_{x=t}^\infty\\
 &=-\frac{p^*}{q}s(\infty)^{-q/p^*}+ \frac{p^*}{q}s(t)^{-q/p^*}\\
 &= \frac{p^*}{q}s(t)^{-q/p^*} \qd\text{if  (\ref{a-07}) holds}.
 \endaligned
 $$
That is,
$$s(t)^{-q/p^*}= \frac{q}{p^*} \int_t^\infty s(x)^{-q/p^*-1}\d s(x)\qd\text{if  (\ref{a-07}) holds}.$$
Combining this with (\ref{a-11}), under (\ref{a-07}), we obtain
$$\aligned
\int_0^\infty \big(H f(x)\big)^q u(x)\d x
&\le \frac{q}{p^*}B^q \int_0^\infty\bigg[\int_t^\infty s(x)^{-q/p^*-1}\d s(x)\bigg]
\d \big(H f(t)\big)^q\\
&=\frac{q}{p^*}B^q \int_0^\infty \bigg[\int_0^x \d \big(H f(t)\big)^q\bigg]
s(x)^{-q/p^*-1}\d s(x)\\
&=\frac{q}{p^*}B^q \int_0^\infty \big(H f(x)\big)^q s(x)^{-q/p^*-1}\d s(x)\\
&=\frac{q}{p^*}B^q \int_0^\infty \big(H \fz(s)\big)^q s^{-q/p^*-1}\d s\qd \text{(by (\ref{a-10}) and (\ref{a-07}))}.
\endaligned$$
Therefore,
$$\bigg[\int_0^\infty \big(H f(x)\big)^q u(x)\d x\bigg]^{1/q}
\le \bigg[\frac{q}{p^*}\bigg]^{1/q}B \bigg[\int_0^\infty \big(H \fz(s)\big)^q s^{-q/p^*-1}\d s\bigg]^{1/q}.$$
By Bliss's Lemma, the right-hand side is controlled by
$$k_{q, p} B \bigg[\int_0^{\infty} \fz(s)^p\d s\bigg]^{1/p}
=k_{q, p} B \bigg[\int_0^{\infty} f(x)^p v(x)\d x\bigg]^{1/p}.$$
We have thus proved the required assertion under (\ref{a-07}).

To remove condition (\ref{a-07}) used in the proof above, we use Lemma \ref{ta-6}.
For given $\mu$ and $\nu$, we define naturally $\mu^N$ and $\nu^N$
to be the restriction of $\mu$ and $\nu$ on $[0, N]$. Then we clearly have
${\tilde\nu}^{N, \nu}=\nu$, respectively. We have already proved that
$$A\big({\tilde\mu}^N, {\tilde\nu}^{N, \d x}\big)\le k_{q,p}
B\big({\tilde\mu}^N, {\tilde\nu}^{N, \d x}\big)$$
since ${\tilde\nu}^{N, \d x}$ satisfies condition (\ref{a-07}).
By Lemma \ref{ta-6}, we get
$$A\big({\tilde\mu}^N, \nu\big)=A_N=A\big({\tilde\mu}^N, {\tilde\nu}^{N, \d x}\big),\qqd
B\big({\tilde\mu}^N, \nu\big)=B_N=B\big({\tilde\mu}^N, {\tilde\nu}^{N, \d x}\big),$$
and then
$$A\big({\tilde\mu}^N, \nu\big)\le k_{q,p}B\big({\tilde\mu}^N, \nu\big).$$
The assertion now follows by letting $N\to\infty$.
\deprf

We conclude the Appendix by a discussion on the eigenequation corresponding to
the Hardy-type inequality,

\prp{\cms Again, let $\mu(\d x)=u(x)\d x$ and $\nu(\d x)=v(x)\d x$.
When $q\ne p$, the eigenequation for the Hardy-type inequality becomes
$$\big(v {g'}^{p-1}\big)'=-u g^{q-1},\qqd g, g'>0\qqd\text{\cms a.e. on } (-M, N)$$
and with boundary condition $\big(v{g'}^{p-1}g\big)\big|_{-M}^{N}=0$ once $M, N<\infty$.
Actually, the eigenequation is equivalent to the following assertion:
$$\frac{u g^{q-1}(x)}{\big(v {g'}^{p-1}\big)'(x)}=:-\ez\qqd
\text{\cms is independent of a.e. $x$ \;$\text{\cms on } (-M, N).$}$$
If the boundary condition holds, then the optimal constant $A$ is given by
$$A=\frac{\|g\|_{L^q(\mu)}}{\|g'\|_{L^p(\nu)}}
=\ez^{1/q} \bigg[\int_{-M}^N v {g'}^p\bigg]^{1/q-1/p}.$$
}\deprp

\prf The first assertion comes from the Euler-Lagrange equation in variational methods. Roughly speaking, the idea goes as follows.
Let $g$ be the minimizer of the inequality and let $h\in {\scr C}_0^\infty(-M, N)$. Define
$$F(\vz)=\bigg(\int_{-M}^N {(g+\vz h)'}^p v\bigg)^{1/p}
\bigg(\int_{-M}^N {(g+\vz h)}^q u\bigg)^{-1/q},\qqd \vz>0.$$
Then, it is easy to check that $\frac{\d}{\d\vz}F(\vz)=0$ iff
$$Q\int_{-M}^N {g'}^{p-1} h' v= P \int_{-M}^N {g}^{q-1}u h, $$
where
$$P=\int_{-M}^N {g'}^p v, \qqd Q=\int_{-M}^N {g}^q u.$$
Using the integration by parts formula, we obtain
$$\int_{-M}^N\bigg[\frac{Q}{P} \big({g'}^{p-1} v\big)'
+{g}^{q-1}u\bigg]h=0.$$
Since $h$ is arbitrary, this gives us
$$\frac{Q}{P} \big({g'}^{p-1} v\big)'
+{g}^{q-1}u=0,\qqd \text{a.e.}$$
Here a key is the inhomogenous, one may replace $g$ by $\xi g$ if necessary for some constant $\xi$, so that the coefficient $Q/P$ can be set
to be one. This gives us the first assertion and than one leads to the equivalent assertion.

Multiplying by $g$ on both sides of the eigenequation, and using the integral by parts formula,
it follows that
$$\big(v {g'}^{p-1}g\big)\big|_{-M}^{N}
-\int_{-M}^N v {g'}^p
=-\ez^{-1} \int_{-M}^N u g^q.$$
By boundary condition, we obtain
$$\int_{-M}^N u g^q=\ez \int_{-M}^N v {g'}^p.$$
Hence
$$\|g\|_{L^q(\mu)}= \ez^{1/q}\|g'\|_{L^p(\nu)}^{p/q} $$
which is the last assertion.
\deprf

To apply the last result to Example \ref{t2-5}, let
$$g(x) = \frac{\az x}{(1 + \bz x^{\gz})^{1/{\gz}}}>0, \qqd \az, \bz>0, \; \gz=\frac q p -1>0.$$
Then
$$g'(x) = \frac{\az}{(\bz x^{\gz} + 1)^{({\gz} + 1)/{\gz}}  }>0.$$
Clearly, we have $\big(g{g'}^{p-1}\big)(0)=0$. Next, since $g(x)\sim 1$
and $g'(x)\sim x^{-q/p}$ as $x\to\infty$, we also have
$\lim_{x\to\infty}\big(g{g'}^{p-1}\big)(x)=0$. Some computations show that
$$\aligned
\big({g'}^{p-1}\big)'(x)& =-
   \frac{\az^{p-1}\bz (p-1)(\gz +1)x^{\gz-1}}{(1+\bz x^{\gz})^{p+(p-1)/\gz}}\\
   &=  -\frac{\az^{p-1}\bz (p-1)(\gz +1)x^{\gz-1}}{(1+\bz x^{\gz})^{(q-1)/\gz}}\qqd\bigg(p+\frac{p-1}{\gz}=\frac{q-1}{\gz}\bigg),\\
  \frac{g(x)^{q-1}}{x^{{\gz} p + p - {\gz}}}&=
  \frac{\az^{q-1} x^{{\gz} - p (1 + {\gz})+q-1}}{(1 + \bz x^{\gz})^{(q-1)/{\gz}}}\\
  &= \frac{\az^{q-1} x^{{\gz} - 1}}{(1+\bz x^{\gz})^{(q-1)/\gz}} \qqd \big(q=p (1 + {\gz})\big).
\endaligned
$$
Hence
$$ \frac{g(x)^{q-1}}{x^{{\gz} p + p - {\gz}}}\bigg/
\big({g'}^{p-1}\big)'(x)
=-\frac{\az^{q-p}}{\bz (p-1)(\gz +1)}=:- \ez.$$
[The right-hand side is independent of $x$ for all $\bz$, for simplicity, one may simply
set $\bz=1$. Then one can also set $\az=1$ in computing $A$. This observation may simplify the computation blow.] Set $t=\bz x^\gz$, then
$$\aligned
\d t&= \bz \gz x^{\gz-1}\d x,\qqd x=\bigg(\frac{t}{\bz}\bigg)^{1/\gz}.\\
{g'}(x)^p \d x
&=\frac{\az^p}{(1 + \bz x^\gz)^{p (1 + \gz)/\gz}}\d x\\
&=\frac{\az^p x^{1-\gz}}{\bz\gz(1 + \bz x^\gz)^{p (1 + \gz)/\gz}}\d t\\
&=\frac{\az^p}{\bz^{1/\gz}\gz} \frac{t^{1/\gz-1}}{ (1 + t)^{p (1 + \gz)/\gz}}\d t.
\endaligned$$
Because
$$\text{\rm B}(x, y)=\int_0^\infty \frac{t^{x-1}}{(1+t)^{x+y}}\d t,$$
we obtain
$$\int_0^\infty {g'}(x)^p \d x
=\frac{\az^p}{\bz^{1/\gz}\gz}\int_0^\infty \frac{t^{1/\gz-1}}{ (1 + t)^{p (1 + \gz)/\gz}}\d t
=\frac{p\az^p\bz^{p/(p-q)}}{q-p}\text{\rm B}\bigg(\frac{p}{q-p}, \frac{p(q-1)}{q-p}\bigg). $$
Therefore, we arrive at
$$\aligned
&\bigg[\int_0^\infty \!\!\frac{g(x)^q}{x^{\gz p + p - \gz}}\d x\bigg]^{1/q}
\bigg/\bigg[\int_0^\infty {g'}(x)^p \d x\bigg]^{1/p}\\
&\qd =\bigg[\frac{\az^{q-p}}{\bz (p-1)(\gz +1)}\bigg]^{1/q}
\bigg[\frac{p\az^p\bz^{p/(p-q)}}{q-p}\text{\rm B}\bigg(\frac{p}{q-p}, \frac{p(q-1)}{q-p}\bigg)\bigg]^{1/q-1/p}\\
&\qd =\bigg[\frac{p^*}{q}\bigg]^{1/q} \Bigg[\frac{q-p}{p\text{\rm B}\big(\frac{p}{q-p}, \frac{q(p-1)}{q-p}\big) }\Bigg]^{1/p-1/q}
\endaligned$$
which is the exact upper bound given in Example \ref{t2-5} or Lemma \ref{ta-4}.
\medskip

\nnd{\bf Acknowledgments}. {\small
The paper was started when the author visited Taiwan in October, 2012.
The invitation and financial support from Inst. Math. Acad. Sin.
(hosts: Chii-Ruey Hwang and Shuenn-Jyi Sheu) are acknowledged. The author thanks
the very warm hospitality made by C.R. Hwang, S.J. Sheu, Tzuu-Shuh Chiang,
Yun-Shyong Chow and their wives. The author also thanks the following
professors for their warm hospitality and the financial support from their
universities: Lung-Chi Chen at Fu Jen Catholic University,
Hong-Kun Xu, Mong-Na Lo Huang, and Mei-Hui Guo at National Sun Yat-sen University,
Yuh-Jia Lee at National University of Kaohsiung,
Shuenn-Jyi Sheu at National Central University,
Yuan-Chung Sheu at National Chiao Tung University,
Tien-Chung Hu at National Tsing Hua University.
In the past years, the author has obtained a lot of help from Y.S. Chow,
as well as the librarians at the Institute of Mathematics
for providing some old papers. Without their help, this paper would not exist.

Research is supported in part by the
         National Natural Science Foundation of China (No. 11131003),
         the ``985'' project from the Ministry of Education in China,
and the Project Funded by the Priority Academic Program Development of Jiangsu Higher Education Institutions.

The careful corrections to an earlier version of the paper by an unknown referee is also acknowledged.}

\bigskip

\nnd {\small School of Mathematical Sciences, Beijing Normal University,
Laboratory of Mathematics and Complex Systems (Beijing Normal University),
Ministry of Education, Beijing 100875,
    The People's Republic of China.\newline E-mail: mfchen@bnu.edu.cn\newline Home page:
    http://math.bnu.edu.cn/\~{}chenmf/main$\_$eng.htm
}

\end{document}